
\magnification1200
\input amstex.tex
\documentstyle{amsppt}
\nopagenumbers
\hsize=12.5cm
\vsize=18cm
\hoffset=1cm
\voffset=2cm

\footline={\hss{\vbox to 2cm{\vfil\hbox{\rm\folio}}}\hss}

\def\DJ{\leavevmode\setbox0=\hbox{D}\kern0pt\rlap
{\kern.04em\raise.188\ht0\hbox{-}}D}

\def\txt#1{{\textstyle{#1}}}
\baselineskip=13pt
\def\hf{{\textstyle{1\over2}}}
\def\a{\alpha}\def\b{\beta}
\def\d{{\,\roman d}}
\def\e{\varepsilon}\def\E{{\roman e}}
\def\f{\varphi}
\def\b{\beta} \def\g{\gamma}
\def\G{\Gamma}

\def\s{\sigma}

\def\={\;=\;}

\def\zt{\zeta(\hf+it)}

\def\D{\Delta}
\def\E{{\roman e}}

\def\z{\zeta}

\def\hf{{\textstyle{1\over2}}}
\def\txt#1{{\textstyle{#1}}}
\def\f{\varphi}

\def\le{\leqslant} \def\ge{\geqslant}
\font\tenmsb=msbm10
\font\sevenmsb=msbm7
\font\fivemsb=msbm5
\newfam\msbfam
\textfont\msbfam=\tenmsb
\scriptfont\msbfam=\sevenmsb
\scriptscriptfont\msbfam=\fivemsb
\def\Bbb#1{{\fam\msbfam #1}}

\def \NN {\Bbb N}

\def \ZZ {\Bbb Z}

\font\ff=cmr8
\def\txt#1{{\textstyle{#1}}}
\baselineskip=13pt

\font\teneufm=eufm10
\font\seveneufm=eufm7
\font\fiveeufm=eufm5
\newfam\eufmfam
\textfont\eufmfam=\teneufm
\scriptfont\eufmfam=\seveneufm
\scriptscriptfont\eufmfam=\fiveeufm
\def\mathfrak#1{{\fam\eufmfam\relax#1}}

\font\tenmsb=msbm10
\font\sevenmsb=msbm7
\font\fivemsb=msbm5
\newfam\msbfam
     \textfont\msbfam=\tenmsb
      \scriptfont\msbfam=\sevenmsb
      \scriptscriptfont\msbfam=\fivemsb
\def\Bbb#1{{\fam\msbfam #1}}

\def \NN {\Bbb N}

\def \ZZ {\Bbb Z}

  \def\rightheadline{{\hfil{\ff On some integrals involving the  divisor problem
 }\hfil\tenrm\folio}}

  \def\leftheadline{{\tenrm\folio\hfil{\ff
   A. Ivi\'c and W. Zhai}\hfil}}
  \def\emptyheadline{\hfil}
  \headline{\ifnum\pageno=1 \emptyheadline\else
  \ifodd\pageno \rightheadline \else \leftheadline\fi\fi}

\topmatter

\title
On certain integrals involving the Dirichlet divisor problem
\endtitle
\author   Aleksandar Ivi\'c and Wenguang Zhai
 \endauthor

\nopagenumbers

\medskip

\address
Aleksandar Ivi\'c, Serbian Academy of Sciences and Arts,
Knez Mihailova 35, 11000 Beograd, Serbia.
\endaddress

\address
Wenguang Zhai, Department of Mathematics, China University of Mining and Technology,
Beijing  100083, China.
\endaddress

\bigskip

\bigskip
\email {\medskip
\tt
aleksandar.ivic\@rgf.bg.ac.rs, aivic\_2000\@yahoo.com;   zhaiwg\@hotmail.com}
\endemail
\dedicatory
\enddedicatory
\abstract
{We prove that
$$
\int_1^X\D(x)\D_3(x)\d x \ll X^{13/9}\log^{10/3}X, \quad \int_1^X\D(x)\D_4(x)\d x \ll_\e X^{25/16+\e},
$$
where $\D_k(x)$ is the error term in the asymptotic formula for the summatory function of $d_k(n)$,
generated by $\z^k(s)$ ($\D_2(x) \equiv \D(x)$). These bounds are sharper than the ones which follow
by the Cauchy-Schwarz inequality and mean square results for $\D_k(x)$. We also obtain the analogues of
the above bounds when $\D(x)$ is replaced by $E(x)$, the error term in the mean square formula for $|\zt|$.
 }
\endabstract
\endtopmatter

\document

\head
1. Introduction and statement of results
\endhead

Let $d_k(n)$ denote the number of ways $n$ can be written as a product of $k$ factors, where $k\in\NN$
is given. Thus
$d_k(n)$ is generated by $\z^k(s)\;(\Re s > 1)$, where $\z(s)$ denotes the familiar zeta-function
of Riemann. It is seen that $d_1(n) \equiv 1$ and $d_2(n) \equiv d(n)$ is the number of divisors of $n$.
 The error term in the asymptotic formula for sums of
$d_k(n)$ is commonly denoted by $\D_k(x)$, namely
$$
\D_k(x) = \sum_{n\le x}d_k(n) - xP_{k-1}(\log x),\leqno(1.1)
$$
where $P_{k-1}(t)$ is a suitable polynomial of degree $k-1$ in $t$, whose coefficients depend on $k$.
Usually one writes $\D_2(x) \equiv \D(x)$, and the estimation of this function is called the Dirichlet
divisor problem, while the estimation of $\D_k(x)$ for $k\ge3$ is known as the general (or generalized)
divisor problem,
or the Piltz divisor problem. A comprehensive account on $\D_k(x)$ is to be found in Chapter 12 of
E.C. Titchmarsh \cite{Tit} or Chapter 13 of A. Ivi\'c \cite{Iv1}. The coefficients of $P_{k-1}(t)$ may be evaluated
by using the formula
$$
P_{k-1}(\log x) = \mathop{\roman{Res}}_{s=1}\, x^{s-1}\z^k(s)s^{-1},
$$
and the fact that
$$
\z(s) = \frac{1}{s-1} + \g + \sum_{k=1}^\infty \g_k(s-1)^k
$$
in a neighborhood of $s=1$, where $\g = -\G'(1) = 0.577215\ldots$ is Euler's constant. In this way one
finds that
$$
\eqalign{
P_1(t) &= t + (2\g-1),\cr
P_2(t) & = \hf t^2+ (3\g-1)t + (3\g^2-3\g+ 3\g_1+1),\cr}
$$
and so on.

\medskip
Following standard notation one defines $\a_k$ and $\b_k$ as the infima of numbers $a_k$ and $b_k$,
respectively, for which
$$
\D_k(x) \,\ll\, x^{a_k},\qquad \int_1^x\D_k^2(y)\d y \,\ll\, x^{1+2b_k}.\leqno(1.2)
$$
It is obvious that $\b_k \le \a_k\;(\forall k )$. It is known that $\b_k\ge (k-1)/(2k)$, and
it is conjectured that
$$
\a_k \= \b_k \= \frac{k-1}{2k}\quad(k>1).\leqno(1.3)
$$
This is not yet proved for any $\a_k$, and for $\b_k$ (1.3) is known to hold only for $k =2,3,4$.
Currently the best known upper bounds for $\a_2, \a_3$ and $\a_4$ are
$$
\a_2 \le \frac{517}{1648} = 0.31371\ldots,\quad \a_3 \le \frac{43}{96} = 0.447916\ldots,\quad\a_4 \le \frac12.
$$
The bound for $\a_2$ is a recent result of J. Bourgain and N. Watt \cite{BW}, the bound for $\a_3$ was proved
in 1981 by G. Kolesnik \cite{Kol}, and the bound for $\a_4$ follows easily by the Perron inversion formula.

It is not only that $\b_2 = 1/2, \b_3 = 1/3$ is known, but actually one has the asymptotic formulas
$$\int_1^X\D^2(x)\d x = A_2X^{3/2} + O(X\log^3X\log\log X)\;\;\Bigl(A_2 =
(6\pi^2)^{-1}\sum_{n=1}^\infty d^2(n)n^{-3/2}\Bigr),
\leqno(1.4)
$$
and
$$
\int_1^X\D_3^2(x)\d x = A_3X^{5/3} + O_\e(X^{14/9+\e})\;\;\Bigl(A_3
= (10\pi^2)^{-1}\sum_{n=1}^\infty d_3^2(n)n^{-4/3}\Bigr).
\leqno(1.5)
$$
The asymptotic formula (1.4) is due to Y.-K. Lau and K.-M. Tsang \cite{LT},
while (1.5) was proved in 1956 by K.-C. Tong
\cite{Ton}. On the other hand, $\D_k(x)$ takes large positive and small negative values.
Indeed, Theorem 2 of A. Ivi\'c \cite{Iv2}
 states that, for $k\ge2$ fixed, there exist constants $B, C>0$ such that for $T\ge T_0$ the interval
$[T, T+ CT^{(k-1)/k}]$ contains two points $t_1, t_2$ for which
$$
\D_k(t_1) > Bt_1^{(k-1)/(2k)},\quad \D_k(t_2) < -Bt_2^{(k-1)/(2k)}.
$$
Therefore it seems reasonable to expect a lot of cancellation in the integral of
$\D_k(x)\D_\ell(x)\;(k\ne \ell; \,k,\ell >1)$.
In other words, the integral of this function should be of a smaller order of magnitude than what one expects
is the order of the integral of $|\D_k(x)| |\D_\ell(x)|$. We are able to establish two results in this direction
when $k=2$. We shall prove

\medskip
THEOREM 1. {\it We have}
$$
\int_1^X\D(x)\D_3(x)\d x \;\ll\; X^{13/9}(\log X)^{10/3}.\leqno(1.6)
$$

\medskip
THEOREM 2. {\it We have}
$$
\int_1^X\D(x)\D_4(x)\d x \;\ll_\e\; X^{25/16+\e}.\leqno(1.7)
$$

\medskip
Here and later $\e$ denotes arbitrarily small positive constants, not necessarily the same ones
at each occurrence. The symbol $a \ll_\e b$ means that the implied $\ll$-constant depends only on $\e$.
To assess the strength of the bounds in (1.6) and (1.7) note that, using (1.4),
(1.5) and the Cauchy-Schwarz inequality for integrals, it follows that
$$
\left|\int_1^X\D(x)\D_3(x)\d x\right| \le \left(\int_1^X\D^2(x)\d x \int_1^X\D_3^2(x)\d x\right)^{1/2}
\ll X^{19/12}.
$$
Similarly, from (1.4) and $\b_4 = 3/8$ we have
$$
\int_1^X\D(x)\D_4(x)\d x \ll_\e X^{13/8+\e}.
$$
Since
$$
13/9 = 1.\bar{4} < 19/12 = 1.58\bar{3}, \quad 25/16 = 1.5625 < 13/8 = 1.625,
$$
we show that there is a substantial cancellation in the integrals in (1.6) and (1.7). However, the
true order of these integrals remains elusive and represents a difficult problem.
We note that the estimation of the integral of $\D(x)\D_\ell(x)$ when $\ell\ge5$ seems difficult,
and that $\b_\ell = (\ell-1)/(2\ell)$ is not known yet for any $\ell\ge5$; for $\ell =5$ the
sharpest bound is $\b_5 \le 9/20$, due to   W.-P. Zhang \cite{Zh}, while it is conjectured that
$\b_5 = 2/5$.

We note that in \cite{IvZh} it was proved that
$$
\int_0^T \D(t)|\zt|^2\d t \;\ll\; T(\log T)^4,\leqno(1.8)
$$
but we remarked that obtaining an asymptotic formula for the integral in (1.8) seems difficult.

\medskip
The method of proof of Theorem 1 is of a fairly general nature,
and works if $\D(x)$ is replaced by a number-theoretic error term which, broadly
speaking, has a structure similar to $\D(x)$. Perhaps the most interesting case is when
$\D(x)$ is replaced by
$$
E(x) := \int_0^x|\zt|^2\d t - x\Bigl(\log\bigl(\frac{x}{2\pi}\bigr) + 2\g -1\Bigr),\leqno(1.9)
$$
the error term in the mean square formula for $|\zt|$. This is a fundamental function in the
theory of $\z(s)$, and the reader is referred to Chapter 15 of \cite{Iv1} and Chapter 3 of \cite{Iv3}
for an account of $E(x)$. We shall prove

\medskip
THEOREM 3. {\it We have}
$$
\int_1^X E(x)\D_3(x)\d x \ll X^{3/2}\log^{3/2}X,\quad
\int_1^X E(x)\D_4(x)\d x \ll_\e X^{25/16+\e}.\leqno(1.10)
$$

\medskip
There are also several other ways to generalize Theorem 1 and Theorem 2, namely by replacing $\D(x)$ with
another suitable number-theoretic error term. These possibilities were already analyzed
in detail in \cite{Iv6}, where the functions $P(x)$ and $A(x)$ were mentioned. One has
$$
P(x) \,:=\, \sum_{n\le x}r(n) - \pi x, \qquad r(n) \,=\, \sum_{n=a^2+b^2}1,
$$
and
$$
A(x)  \;:=\; \sum_{n\le x}a(n),
$$
where $a(n)$ is the $n$-th Fourier
coefficient of $\varphi(z)$, and $\varphi(z)$ is a
holomorphic cusp form of weight $\kappa$ with respect to the full
modular group $SL(2,\ZZ)$. However, correlating these functions with $\D_k(x)$ seems a little
artificial, while the integrals of $\D(x)\D_k(x)$ seemed much more natural to investigate.

\smallskip
{\bf Acknowledgement.} W. Zhai is  supported by the National Key Basic Research
Program of China (Grant No. 2013CB834201).

\bigskip
\head
2. The necessary lemmas
\endhead
\medskip
In this section we formulate the lemmas needed for the proof of our results. The first lemma gives
a good expression for $\D(x)$ and a mean square estimate involving the ``tails'' of this expression.
It is the lack of such a result for $\D_k(x)$ when $k\ge3$ that thwarts the efforts to obtain
satisfactory bounds for the integral of $\D_k(x)\D_\ell(x)$ when $3 \le k < \ell$.

\medskip
{\bf Lemma 1}. We have, for $1 \ll N \ll x$,
$$
\D(x) = {1\over\pi\sqrt{2}}x^{1\over4}
\sum_{n\le N}d(n)n^{-{3\over4}}
\cos(4\pi\sqrt{nx} - {\txt{1\over4}}\pi) + R(x), \leqno(2.1)
$$
where
$$
R(x) \,=\, R(x,N) \,=\, O_\e(x^{{1\over2}+\e}N^{-{1\over2}}).\leqno(2.2)
$$
 If $1\ll N\ll X\log^{-4} X, $ then
$$
\int_X^{2X}R^2(x,N)\d x \;\ll\; X^{3/2}N^{-1/2}\log^3 X.\leqno(2.3)
$$

\medskip
{\bf Proof}. The formulas (2.1)--(2.2) are the classical truncated Vorono{\"\i} formula for $\D(x)$.
For this see e.g., Chapter 3 of A. Ivi\'c \cite{Iv1} for a proof. To obtain (2.3), note that by a result
of T. Meurman \cite{Meu} one has, for $Q \gg x \gg 1$,
$$
\D(x) = {1\over\pi\sqrt{2}}x^{1\over4}
\sum_{n\le Q}d(n)n^{-{3\over4}}
\cos(4\pi\sqrt{nx} - {\txt{1\over4}}\pi) + F(x),
\leqno(2.4)
$$
where $F(x) \ll x^{-1/4}$ if $||x|| \gg x^{5/2}Q^{-1/2}$, and we always have $F(x) \ll_\e x^\e$.
Here, as usual, $||x||$ is the distance of the real number $x$ to the nearest integer.

In (2.4) we take $Q=X^7, X\le x\le 2X$. Then, for $1 \ll N \ll X$, $R(x)$ in (2.1) can be written as
$$
R(x) = {1\over\pi\sqrt{2}}x^{1\over4}
\sum_{N<n\le Q}d(n)n^{-{3\over4}}
\cos(4\pi\sqrt{nx} - {\txt{1\over4}}\pi) + F(x),
$$
where in our case (2.4) shows that $F(x) \ll_\e x^\e$ always, and $F(x) \ll x^{-1/4}$ if
$||x|| \gg X^{-1}$. Then the integral in (2.3) is
$$
\ll X^{1/2}\int_X^{2X}\Bigl|\sum_{N<n\le Q}d(n)n^{-3/4}\E^{4\pi i\sqrt{nx}}\Bigr|^2\d x
+ \int_X^{2X} F^2(x)\d x.
$$
We have
$$
\eqalign{
\int_X^{2X} F^2(x)\d x & = \int_{X,||x||\gg X^{-1}}^{2X}F^2(x)\d x
+ \int_{X,||x||\ll X^{-1}}^{2X}F^2(x)\d x\cr&
\ll \int_X^{2X}x^{-1/2}\d x + X^{2\e} \ll X^{1/2}.\cr}
$$
Now  squaring out and integrating it follows that
$$
\eqalign{&
\int_X^{2X} \Bigl|\sum_{N<n\le Q}d(n)n^{-3/4}\E^{4\pi i\sqrt{nx}}\Bigr|^2\d x\cr&
\ll X\sum_{n>N}d^2(n)n^{-3/2} + \sum_{N<m\ne n\le Q}\frac{d(m)d(n)\sqrt{X}}{(mn)^{3/4}|\sqrt{m}-\sqrt{n}|}.
\cr}
$$

Note that $\sum\limits_{n\le y}d^2(n)\ll y\log^3 y$ for $y\ge 3.$
Thus  partial summation gives
$$X\sum_{n>N}d^2(n)n^{-3/2}\ll XN^{-1/2}\log^3 X.$$
The second sum above does not exceed $\,S_1+S_2,$ where
$$S_j=\sum_{j}\frac{d(m)d(n)\sqrt{X}}{(mn)^{3/4}|\sqrt{m}-\sqrt{n}|}\qquad (j=1,2),
$$
and  (SC($A$) means: summation conditions for $A$)
$$
\eqalign{&
\roman{SC}\,(S_1): N<m\ne n\le Q,  |\sqrt{m}-\sqrt{n}|\le (mn)^{1/4}/10,\cr&
\roman{SC}\,(S_2): N<m\ne n\le Q,  |\sqrt{m}-\sqrt{n}|\ge (mn)^{1/4}/10.
\cr}
$$
The sum $S_2$ is bounded by
$$ \sqrt{X}\sum_{m,n\le Q}\frac{d(m)d(n)}{mn}\ll \sqrt{X}\log^4 X$$
by using the well-known estimate $\sum\limits_{n\le y}d(n)\ll y\log y.$
In $S_1$ we have
$m\asymp n.$ By
$d(m)d(n) \le \hf(d^2(m)+d^2(n))$  we have
$$
\eqalign{
S_2&\ll X^{1/2} \sum_{N<m\le Q}\frac{d^2(m)}{m}\sum_{N<n\le Q,n\ne m}\frac{1}{|n-m|}\cr&
\ll X^{1/2} \log X\sum_{N<n\le Q}\frac{d^2(m)}{m} \ll X^{1/2}\log^5X.\cr}
$$
Collecting the above estimates we get  the bound in (2.3), since
$$
X\log^5X \ll X^{3/2}N^{-1/2}\log^3X\qquad(1\ll  N \ll X\log^{-4}X).
$$

\medskip
{\bf Lemma 2}. Let $1/6 \le \s_0 < 1/2\,$ be a constant. Then, for $k = 3$ and $k = 4$,  we have
$$
\D_k(x) = \frac{1}{2\pi i}\int_{\s_0-iT}^{\s_0+iT}\frac{\z^k(s)}{s}x^s\d s + O_\e(X^{1+\e}T^{-1})
+ O(x^{\s_0}T^{k(1/2-\s_0)-1}),\leqno(2.5)
$$
where $X \le x \le 2X, 1 \ll T \ll X$.

\medskip
{\bf Proof}. From the Perron inversion formula (see e.g., (A.10) of \cite{Iv1}) we have, for $X\le x \le 2X,
1\ll T \ll X$ and sufficiently small $\e>0$,
$$
\D_k(x) = \frac{1}{2\pi i}\int_{1-\e-iT}^{1-\e+iT}\frac{\z^k(s)}{s}x^s\d s + O_\e(X^{1+\e}T^{-1}).
$$
By Cauchy's theorem we replace the segment of integration $[1-\e-iT, 1-\e+iT]$
by the segment $[\s_0-iT, \s_0+iT]$.
In this process we are making an error which is $\ll J_1(T) + J_2(T)$, where
$$
J_1(T) := \int_{\s_0}^{1/2}\frac{x^\s}{T}|\z(\s+iT)|^k\d\s, \;
J_2(T) := \int_{1/2}^{1}\frac{x^\s}{T}|\z(\s+iT)|^k\d\s.
$$
Then
$$
\eqalign{
\int_T^{2T}J_1(t)\d t & = \int_{\s_0}^{1/2}x^\s\left(\int_T^{2T}\frac{|\z(\s+it)|^k}{t}\d t\right)\d\s
\cr&
\ll \int_{\s_0}^{1/2}x^\s T^{k(1/2-\s)}\left(\int_T^{2T}\frac{|\z(1-\s+it)|^k}{t}\d t\right)\d\s
\cr&
\ll \int_{\s_0}^{1/2}x^\s T^{k(1/2-\s)}\d\s \ll X^{1/2} + X^{\s_0}T^{k(1/2-\s_0)}.
\cr}\leqno(2.6)
$$
Here we used the functional equation
$$
\z(s) = \chi(s)\z(1-s), \quad |\chi(s)| =\Bigl|\frac{\G(\hf(1-s))}{\G(\hf s)}\pi^{s-1/2}\Bigr| \asymp |t|^{1/2-\s},
\leqno(2.7)
$$
and the bound (follows from Lemma 5)
$$
\int_T^{2T}|\z(\a+it)|^k\d t \,\ll\, T\qquad(1/2 < \a \le 1; \,k = 3,4).\leqno(2.8)
$$
Note that (2.8) is not known to hold in the whole range $1/2 < \a \le 1$ when $\a$ is a constant,
and $k > 4$ is an integer, which is one of the reasons why obtaining the analogue of Theorem 1 and 2
for the integral of $\D(x)\D_k(x)$ when $k >4$ is difficult.

Similarly to (2.6) we obtain
$$
\eqalign{&
\int_T^{2T}J_2(t)\d t \ll \int_{1/2}^{1}x^\s\left(\int_T^{2T}\frac{|\z(\s+it)|^k}{t}\d t\right)\d\s\cr&
\ll \int_{1/2}^{1}x^\s\d\s \ll X.\cr}
$$
It follows that
$$
\int_T^{2T}(J_1(t)+ J_2(t))\d t \ll X + X^{\s_0}T^{k(1/2-\s_0)}.
$$
This means that there exists $T_0 \in [T, 2T]$ such that
$$
J_1(T_0)+ J_2(T_0)\ll X/T_0 + X^{\s_0}T_0^{k(1/2-\s_0)-1}.
$$
If instead of the initial $T$ we take this $T_0$ and call it again $T$, we obtain (2.5).

\medskip
{\bf Lemma 3}. Let $f(x), \f(x)$ be real-valued functions on $[a,b]$ which satisfy

1. $f^{(4)}(x)$ and $\f^{(2)}(x)$ are continuous;

2. there exist numbers $H, U, A, 0 < H, A < U, 0 < b-a \le U$ such that
$$
\eqalign{&
A^{-1} \ll f^{(2)}(x) \ll A^{-1}, f^{(3)}(x) \ll A^{-1}U^{-1}, f^{(4)}(x) \ll A^{-1}U^{-2},
\cr&
\f(x) \ll H, \f^{(1)}(x) \ll HU^{-1}, \f^{(2)}(x) \ll HU^{-2};
\cr}
$$

3. for some $c$, $a \le c \le b$, $f'(c) = 0$. Then

$$
\eqalign{&
\int_a^b \f(x)\exp\bigl(2\pi if(x)\bigr)\d x  = \frac{1+i}{\sqrt{2}}\cdot \frac{\f(c)\exp(2\pi ic)}{\sqrt{f''(c}}
+ O(HAU^{-1})\cr&
+ O\Bigl(H\min(|f'(a)|^{-1}, \sqrt{A})\Bigr) + O\Bigl(H\min(|f'(b)|^{-1}, \sqrt{A})\Bigr).
\cr}
\leqno(2.9)
$$

\medskip
This is a version of the classical result on exponential integrals with a ``saddle'' point $c$ (see e.g.,
\cite{Iv1} and \cite{Tit}). The particular version embodied in (2.9) is Lemma 2 on p. 71 of
 the monograph of A.A. Karatsuba and S.M. Voronin \cite{KV}. The proof actually shows that there is no main
 term in (2.9) if $c \not \in (a,b)$. If $c = a$ or $c =b$, then the respective main term is to be halved.

\medskip
{\bf Lemma 4}. For $\hf\le\s < 1$ fixed, $1\ll x, y\ll t^k, s = \s+it, xy = \bigl(\frac{t}{2\pi}\bigr)^k$,
$t \ge t_0$ and $k\ge1$ a fixed integer, we have
$$
\eqalign{
\z^k(s) &= \sum_{m=1}^\infty \rho(m/x)d_k(m)m^{-s} + \chi^k(s)\sum_{m=1}^\infty \rho(m/y)d_k(m)m^{s-1}
\cr& + O(t^{k(1-\s)/3-1} + O(t^{k(1/2-\s)-2}y^\s\log^{k-1}t).\cr}\leqno(2.10)
$$
Here $\chi(s)$ is given by (2.7), and  $\rho(u)\;(\ge 0)$ is a smooth function such that $\rho(0)=1, \rho(u) = 0$ for $u \ge 2$.

\medskip
This is Theorem 4.2 of A. Ivi\'c \cite{Iv3}. The explicit construction of $\rho(u)$ is given by Lemma 4.3
therein. The point of smoothing is to have much better error terms than those which can be at present
obtained without it; see e.g., Theorem 4.3 of \cite{Iv1}.

\medskip
{\bf Lemma 5}. For fixed $\s$ such that $\hf < \s \le 1$, we have
$$
\int_1^T|\z(\s+it)|^4\d t \= \frac{\z^2(2\s)}{\z(4\s)}T + O(T^{2-2\s}\log^3T).\leqno(2.11)
$$

\medskip
This is Theorem 2 of A. Ivi\'c \cite{Iv4}.

\bigskip
\head
3. Proof of Theorem 1
\endhead
We may consider, both in the proof of Theorem 1 and Theorem 2, that the integration is over
$[X, 2X]$. Namely if we obtain the bounds in (1.6) and (1.7) for such an integral, then replacing
$X$ by $X2^{-j}$, summing over $j$ and adding all the results
we easily obtain the assertions of (1.6) and (1.7), respectively.

We begin by noting that the contribution of $R(x)$ in (2.1) to
$$
I(X) := \int_X^{2X}\D(x)\D_3(x)\d x \leqno(3.1)
$$
is, by the Cauchy-Schwarz inequality, (1.5) and (2.3),
$$
\ll \left(\int_X^{2X}R^2(x)\d x\int_X^{2X}\D^2_3(x)\d x\right)^{1/2} \ll X^{19/12}N^{-1/4}\log^{3/2}X,
\leqno(3.2)
$$
if $1 \ll N \ll X\log^{-4}X$. By (2.1) and (2.2) there remains a multiple of
$$
\int_X^{2X}x^{1\over4}\sum_{n\le N}d(n)n^{-{3\over4}}\cos(4\pi\sqrt{nx} - {\txt{1\over4}}\pi)\D_3(x)\d x.
\leqno(3.3)
$$
For $\D_3(x)$ we use (2.5) of Lemma 2 with $k=3$. The contribution of the error terms to (3.3) will be
$$
\ll_\e X^{1/4}(X^{1+\e}T^{-1} + X^{\s_0}T^{1/2-3\s_0})\int_X^{2X}\Bigl|\sum_{n\le N}d(n)n^{-{3/4}}
\E^{4\pi i\sqrt{nx}}\Bigr|\d x.
$$

Henceforth we assume that
$$
T\= X, \quad\s_0 \= \frac14,\leqno(3.4)
$$
so that $X^{1+\e}T^{-1} + X^{\s_0}T^{1/2-3\s_0} \ll_\e X^\e$.
By the method of proof of Lemma 1 it is seen that the sum over $n$ in (3.3) is $\ll 1$ in mean square,
if $N\ll X$. This means that the contribution of the error terms to (3.3) is
$$
\ll_\e\;X^{5/4+\e}\leqno(3.5)
$$
if (3.4) holds. Since $z + \bar{z} = 2\Re z$,  we are left with
$$
\eqalign{&
\sum_{n\le N}d(n)n^{-3/4}\int_{-X}^X\frac{\z^3(\frac14 +it)}{\frac14 +it}\int_X^{2X}x^{1/2+it}
\cos(4\pi\sqrt{nx} - {\txt{1\over4}}\pi)\d x\d t\cr&=
2\sum_{n\le N}d(n)n^{-3/4}\Re \left\{\int_0^X\frac{\z^3(\frac14 +it)}{\frac14 +it}\int_X^{2X}x^{1/2+it}
\cos(4\pi\sqrt{nx} - {\txt{1\over4}}\pi)\d x\d t\right\}.\cr}\leqno(3.6)
$$

It transpires that  integration over $x$ leads to exponential integrals of the form
$$
\int_X^{2X}x^{1/2}\E^{iF(x)}\d x,\leqno(3.7)
$$
where henceforth $\,t>0$ and
$$
F(x) \= F(x;t,n) \;:=\; 4\pi\sqrt{nx} \pm t\log x.
$$

We first consider the case of the plus sign. In this case
the derivative
$$
F'(x) = 2\pi\sqrt{n/x} + t/x
$$
is positive. Thus we can apply the first derivative test (Lemma 2.1 of \cite{Iv1}) to the integral in (3.7).
The contribution of (3.7) for $t< \sqrt{nX}$ is $\ll X{n}^{-1/2}$, and its contribution to (3.6) is
$$
\eqalign{&
\ll X\sum_{n\le N}d(n)n^{-5/4}\int_0^{\sqrt{nX}}\frac{t^{3/4}|\z(\frac34 +it)|^3}{|\frac14+it|}\d t\cr&
\ll X\sum_{n\le N}d(n)n^{-5/4}(nX)^{3/8} \ll X^{11/8}N^{1/8}\log X,\cr}
$$
where we used (2.7) and (2.8).
Similarly, for $t > \sqrt{nX}$, one may use again the first derivative test to obtain a contribution which is
$$
\eqalign{&
\ll X^{3/2}\sum_{n\le N}d(n)n^{-3/4}\int_{\sqrt{nX}}^X t^{3/4-1}\frac{|\z(\frac34 +it)|^3}{t}\d t
\cr&
\ll
X^{3/2}\sum_{n\le N}d(n)n^{-3/4}(nX)^{-1/8} \ll  X^{11/8}N^{1/8}\log X.
\cr}
$$

From now on we consider in detail   the case of the minus sign. We have
$$
\eqalign{
F(x) &= 4\pi\sqrt{nx} - t\log x,\cr
F'(x) &= 2\pi n^{1/2}x^{-1/2} - tx^{-1},\cr
F''(x)& = tx^{-2} - \pi n^{1/2}x^{-3/2}.
\cr}
$$

For $t < c_2\sqrt{nX}$ ($c_j$ denotes positive constants) we have $|F'(x)| \gg n^{1/2}X^{-1/2}$, and as in
the previous case the contribution will be
$$
\ll X^{11/8}N^{1/8}\log X.\leqno(3.8)
$$
The bound in (3.8) will also hold for the contribution of $t$ satisfying $t > c_1\sqrt{nX}$, when
$|F'(x)| \gg t/X$. Actually we can take $c_2=2\pi (1-\varepsilon), c_1=2\sqrt{2}\pi (1+\varepsilon).$

The integral in (3.7) may have a saddle point (a solution of $F'(x) = 0$) if $2\pi n^{1/2}x_0^{-1/2} - t/x_0=0$,
hence when
$$
x_0 \= \frac{t^2}{4\pi^2 n}.\leqno(3.9)
$$
It is readily checked that $X\le x_0 \le 2X$ for $2\pi\sqrt{nX} \le t\le2\pi\sqrt{2nX}$. The integral
in (3.7) is evaluated by Lemma 4, where one can take
$$
a = X, b = 2X, U = X, \f(x) = x^{1/2}, H = X^{1/2}, f(x) = \frac{F(x)}{2\pi}, A = X^{3/2}n^{-1/2},
$$
and the conditions on $f^{(3)}(x)$ and  $f^{(4)}(x)$  hold. The total contribution of the error term
$O(HAU^{-1})$ is
$$
\eqalign{&
\ll X\sum_{n\le N}d(n)n^{-5/4}\int_{c_1{\sqrt{nX}}}^{c_2{\sqrt{nX}}} t^{-1/4}|\z(\txt{\frac34} +it)|^3\d t
\cr&
\ll X\sum_{n\le N}d(n)n^{-5/4}(nX)^{3/8} \ll X^{11/8}N^{1/8}\log X,
\cr}
$$
which is similar to (3.8).

We consider now the contribution of the error terms
$$
\eqalign{&
O\left\{X^{1/2}\min\left(\frac{X^{3/4}}{n^{1/4}},{\Bigl|\frac{t}{X}-
\frac{2\pi\sqrt{n}}{\sqrt{X}}\Bigr|}^{-1}\right)\right\} +
\cr&
O\left\{X^{1/2}\min\left(\frac{X^{3/4}}{n^{1/4}},{\Bigl|\frac{t}{2X}-
\frac{2\pi\sqrt{n}}{\sqrt{2X}}\Bigr|}^{-1}\right)\right\}.
\cr}
$$
They are treated analogously, so only the second one will be considered in detail. Let
$$
{\Cal I} \,:=\,\bigl [c_1\sqrt{nX}, \, c_2\sqrt{nX}\,\bigr]
$$
be the remaining interval in the $t$-integral. Let further
$$
{\Cal I}_0 \,:=\, \left\{(t\in {\Cal I})\,\wedge\, \left(\frac{X^{3/4}}{n^{1/4}} \le {\Bigl|\frac{t}{2X}-
\frac{2\pi\sqrt{n}}{\sqrt{2X}}\Bigr|}^{-1}\right)\right\},
$$
and for $j = 1,2,\ldots$
$$
{\Cal I}_j \,:=\, \left\{(t\in {\Cal I})\,\wedge\, \left(\frac{2^{j-1}n^{1/4}}{X^{3/4}} \le \Bigl|\frac{t}{2X}-
\frac{2\pi\sqrt{n}}{\sqrt{2X}}\Bigr| \le \frac{2^{j}n^{1/4}}{X^{3/4}}\right)\right\},
$$
so that $j \ll \log X$.
It is easy to see ${\Cal I} =\cup_{j\ge0}{\Cal I}_j$ and
$$
|{\Cal I}_0| \le (nX)^{1/4},\; |{\Cal I}_j| \le 2^{j+1}(nX)^{1/4}\quad(j\ge1),
$$
where
$|\Cal A|$ denotes the cardinality of the set $\Cal A$.

  The contribution of ${\Cal I}_0$ to the integral
over $t$ in (3.6) is
$$
\eqalign{&
\ll \sum_{n\le N}d(n)n^{-3/4}\int_{{\Cal I}_0}\frac{|\z(\frac{1}{4} +it)|^3}{t}\cdot \frac{X^{5/4}}{n^{1/4}}\d t
\cr&
\ll X^{5/4}\sum_{n\le N}d(n)n^{-1}\int_{{\Cal I}_0}t^{-1/4}|\z(\textstyle{\frac{3}{4}} +it)|^3\d t\cr&
\ll X^{5/4}\sum_{n\le N}d(n)n^{-1}(nX)^{-1/8}\int_{{\Cal I}_0}|\z(\textstyle{\frac{3}{4}} +it)|^3\d t.\cr}
$$
Now we invoke Lemma 5 (with $\s = 3/4$) and use H\"older's inequality for integrals to obtain
$$
\eqalign{&
\int_{{\Cal I}_0}|\z(\txt{\frac{3}{4}} +it)|^3\d t \le
\left(\int_{{\Cal I}_0}|\z(\textstyle{\frac{3}{4}} +it)|^4\d t\right)^{3/4}|{\Cal I}_0|^{1/4}
\cr&
\ll \left((nX)^{1/4} + (nX)^{1/4} \log^3 X\right)^{3/4}(nX)^{1/16} \ll (nX)^{1/4}(\log X)^{9/4}.
\cr}
$$
Thus we finally see that the contribution is
$$
\ll X^{11/8}\sum_{n\le N}d(n)n^{-7/8}(\log X)^{9/4} \ll X^{11/8}N^{1/8}(\log X)^{13/4}.
$$
In a similar vein it is shown that the above bound also holds for the contribution of $\cup_{j\ge 1}{\Cal I}_j$,
only with an additional log-factor since $j \ll \log X$.

\medskip
We turn now to the contribution of the saddle points. By Lemma 3 the main contribution is a multiple of
$$
x_0^{1/2}|F''(x_0)|^{-1/2}\E^{iF(x_0)}.
$$
We have, since $x_0 = t^2/(4\pi^2 n),$
$$
F''(x_0) = - \pi n^{1/2}t^{-3}\cdot8\pi^3n^{3/2} + t^{-3}\cdot 16\pi^4n^2 = 8\pi^4t^{-3}n^2.
$$
Thus
$$
|F''(x_0)|^{-1/2} = (8\pi^4)^{-1/2} t^{3/2}n^{-1} \asymp X^{3/4}n^{-1/4},
$$
since $t \asymp \sqrt{nX}$ ($a\asymp b$ means that $a \ll b \ll a$). We also have
$$
F(x_0) = 4\pi\sqrt{\frac{nt^2}{4\pi^2n}}-t\log\left(\frac{t^2}{4\pi^2n}\right) = 2t - 2t\log t + t\log(4\pi^2n).\leqno(3.10)
$$
With this in mind, we are left with a multiple of
$$
\sum_{n\le N}d(n)n^{-3/4}\int_{2\pi\sqrt{nX}}^{2\pi\sqrt{2nX}}x_0^{1/2}\,\frac{\E^{iF(x_0)}}{\sqrt{|F''(x_0)|}}
\cdot\frac{\z^3(\frac14+it)}{\frac14 +it}\d t,\leqno(3.11)
$$
since $x_0 \in (X, 2X)$ exactly for $2\pi\sqrt{nX} < t < 2\pi\sqrt{2nX}$.

\medskip
At this point we use the approximate functional equation for $\z^3(s)$, writing first $\z(s) = \chi(s)\z(1-s)$
with $s = 1/4+it$, noting that $\chi(s)\chi(1-s) \equiv 1$ and that in Lemma 4 one has $t = \Im s>0$.
Thus by (2.10) with $k=3$ we obtain
$$
\eqalign{&
\z^3(s) = \chi^3(s)\z^3(1-s) = \chi^3(s)\overline{\z^3(1-{\bar s})}\cr&
= \chi^3(s)\sum_{m=1}^\infty\rho(m/x)d_3(m)m^{s-1} + \sum_{m=1}^\infty\rho(m/y)d_3(m)m^{-s}\cr&
+ O(1) + O(t^{-2}y^{3/4}\log^2X),\cr}\leqno(3.12)
$$
where $\rho(u)\; (\ge 0)$ is a smooth function such that $\rho(0)=1, \rho(u)=0$ for $u \ge 2$, and
$$
xy = {\left(\frac{t}{2\pi}\right)}^3\qquad(x \gg 1, y \gg 1).\leqno(3.13)
$$
We choose $x = y = (t/(2\pi))^{3/2}$ in (3.13). Then $t^{-2}y^{3/4}\log^2X \ll 1$, so that the
contribution of the error terms in (3.12) to (3.11) will be
$$
\eqalign{&
\ll X^{1/2}\sum_{n\le N}d(n)n^{-3/4}X^{3/4}n^{-1/4}\int_{2\pi\sqrt{nX}}^{2\pi\sqrt{2nX}}\frac{\d t}{t}
\cr& \ll
X^{5/4}\log^2 X.\cr}
$$
There remains
$$
\sum_{n\le N}d(n)n^{-3/4}(I'+I''),
\leqno(3.14)
$$
say, where we put ($s =1/4+it)$
$$
\eqalign{
I'&:= \int_{2\pi\sqrt{nX}}^{2\pi\sqrt{2nX}}x_0^{1/2}\E^{iF(x_0)}|F''(x_0)|^{-1/2}\chi^3(s)s^{-1}\sum_{m\le 2x}
\rho(m/x)d_3(m)m^{s-1}\d t,\cr
I''&:= \int_{2\pi\sqrt{nX}}^{2\pi\sqrt{2nX}}x_0^{1/2}\E^{iF(x_0)}|F''(x_0)|^{-1/2}s^{-1}\sum_{m\le 2x}
\rho(m/x)d_3(m)m^{-s}\d t.\cr}
$$
Since, by (2.7), $\chi(s)$ is essentially a quotient of two gamma-factors, it admits by Stirling's
formula a full asymptotic expansion in term of negative powers of $t$. Thus the main contribution coming
from $\chi^3(s)$ with $s = 1/4+it$ will be (see e.g., (1.9) of [Iv3])
$$
\left(\frac{2\pi}{t}\right)^{3(-1/4+it)}\E^{3i(t+\pi/4)}.
$$
Consequently the dominating terms in $I'$ and $I''$ will be a multiple of (we set
$M = \max(2\pi\sqrt{nX},2\pi m^{2/3}2^{-2/3})$ for shortness)
$$
\sum_{m\le 2 (nX)^{3/4}}d_3(m)m^{-3/4}\int_M^{2\pi\sqrt{2nX}}x_0^{1/2}\left(\frac{t}{2\pi}\right)^{3/4}
\rho\left(\frac{m}{x}\right)\,\frac{\E^{iH_1(t)}}{s\sqrt{|F'' (x_0)|}}\,\d t\leqno(3.15)
$$
and
$$
\sum_{m\le 2 (nX)^{3/4}}d_3(m)m^{-1/4}\int_M^{2\pi\sqrt{2nX}}x_0^{1/2}
\rho\left(\frac{m}{x}\right)\,\frac{\E^{iH_2(t)}}{s\sqrt{|F'' (x_0)|}}\,\d t,\leqno(3.16)
$$
respectively, where
$$
\eqalign{
H_1(t) &:= 3t\log\left(\frac{2\pi \E}{t}\right) + t\log m + F(x_0),\cr
H_2(t) &:= -t\log m + F(x_0).\cr}
$$
We see that (recall (3.10))
$$
\eqalign{
H_1'(t)&= -5\log t - 3 + 3\log(2\pi \E) +\log(4\pi^2mn)=-5\log t+\log (2\pi)^5mn,\cr
H_2'(t)& = -2\log t +  \log\left(\frac{4\pi^2n}{m}\right),\cr}
$$
so that, in our range of $m,n$ and $t$,
$$
H_1'(t) \,\gg\,\log t,\quad H_2'(t) \,\gg\,\log t.
$$
Therefore we can estimate the integrals in (3.15) and (3.16)  by the first derivative test to obtain
$$
\eqalign{
I' &\ll \sum_{m\le 2(nX)^{3/4}}d_3(m)m^{-3/4}X^{1/2}X^{3/4}n^{-1/4}(nX)^{-1/8}(\log X)^{-1}\cr&
\ll (nX)^{3/16}X^{9/8}n^{-3/8}\log X = X^{21/16}n^{-3/16}\log X,\cr}
$$
and
$$
\eqalign{
I'' &\ll \sum_{m\le 2(nX)^{3/4}}d_3(m)m^{-1/4}X^{1/2}X^{3/4}n^{-1/4}(nX)^{-1/2}(\log X)^{-1}\cr&
\ll (nX)^{9/16}X^{3/4}n^{-3/4}\log X = X^{21/16}n^{-3/16}\log X.\cr}
$$
This gives for the expression in (3.14)
$$
\sum_{n\le N}d(n)n^{-3/4}(I'+ I'') \ll X^{21/16}N^{1/16}\log^2X.
$$
Collecting all the estimates we have
$$
I(X) \ll X^{19/12}N^{-1/4}\log^{3/2}X + X^{11/8}N^{1/8}\log^{17/4} X + X^{21/16}N^{1/16}\log^2X.
$$
So with the choice
$$
N \= X^{5/9}(\log X)^{-22/3}
$$
it follows that
$$
I(X) \,\ll\, X^{13/9}(\log X)^{10/3},
$$
which finishes the proof of Theorem 1.

\head
4. Proof of Theorem 2
\endhead
The proof of Theorem 2 is on the same lines as the proof of Theorem 1, so we shall
indicate only the salient points in the proof. This seemed preferable than considering
the general integral of $\D(x)\D_k(x)$, and then distinguishing between the cases $k=3$
and $k=4$. Whenever possible we shall retain the same notation as in the proof of Theorem 1.

\medskip
We recall first that W.-P. Zhang \cite{Zh} proved that $\b_4 = 3/8$, in other words that
$$
\int_X^{2X}\D_4^2(x)\d x \;\ll_\e\; X^{7/4+\e}.\leqno(4.1)
$$
An asymptotic formula for the integral in (4.1), analogous to (1.4) and (1.5), is not known to hold yet.
From (1.4), (4.1) and the Cauchy-Schwarz inequality we have
$$
I(X) := \int_X^{2X}\D(x)\D_4(x)\d x \ll_\e X^{13/8+\e},\leqno(4.2)
$$
and similarly
$$
\int_X^{2X}\D_4(x)R(x)\d x \ll_\e X^{13/8+\e}N^{-1/4},\leqno(4.3)
$$
where $R(x)$ is as in (2.1)--(2.2). By using (2.1) it is seen that there remains a multiple of
$$
\int_X^{2X}x^{1/4}\sum_{n\le N}d(n)n^{-3/4}\cos(4\pi\sqrt{nx}-\pi/4)\D_4(x)\d x.\leqno(4.4)
$$
For $\D_4(x)$ we use Lemma 2 with $k=4, \s_0 = 1/4, T = X$. The error terms in (2.5) are
$\ll_\e X^{1/4}\log^4X$, and their contribution to (4.4) is
$$
\ll X^{3/2}\log^5X.\leqno(4.5)
$$
Thus we are left with
$$
\sum_{n\le N}d(n)n^{-3/4}\int_{-X}^X\frac{\z^4(\frac14+it)}{\frac14+it}\int_X^{2X}x^{1/2+it}
\cos(4\pi\sqrt{nx}-\pi/4) \d x\,\d t.\leqno(4.6)
$$
We are  led again to exponential integrals of the form
$$
\int_X^{2X}x^{1/2}\E^{iF(x)}\d x,\; F(x) = 4\pi\sqrt{nx} \pm t\log x.
$$
As in the case of Theorem 1, one should take care only of the case of the minus sign in $F(x)$. Namely,
in case of the plus sign we use the first derivative test and see that the total contribution is
$$
\ll X^{3/2}N^{1/4}\log X, \leqno(4.7)
$$
and the relevant contribution of the $t$-integral is for $t\asymp \sqrt{nX}$, as in Theorem 1.
The saddle point is again
$$
x_0\= \frac{t^2}{4\pi^2n}.
$$
Lemma 3 is used again. The contribution of the error term $O(HA/U)$ is
$$
\eqalign{&
\ll X\sum_{n\le N}d(n)n^{-5/4}\int_{c_1\sqrt{nX}}^{c_2\sqrt{nX}}|\z(\textstyle\frac{3}{4}+it)|^4\d t\cr&
\ll X^{3/2}N^{1/4}\log X,\cr}
$$
like in (4.7). For the other two error terms in (2.9) we again define the sets ${\Cal I}_j\;(j\ge0)$
as in the proof of Theorem 1.
The contribution of ${\Cal I}_0$ is
$$
\eqalign{&
\ll \sum_{n\le N}d(n)n^{-3/4}\int_{{\Cal I}_0}|\z(\textstyle{\frac{3}{4}}+it)|^4X^{1/2}X^{3/4}n^{-1/4}\d t
\cr&
\ll X^{5/4}\sum_{n\le N}d(n)n^{-1}(|{\Cal I}_0| + (nX)^{1/4}\log^3 X)\cr&
\ll X^{3/2}N^{1/4}\log^4 X.
\cr}
$$
This is actually simpler than the analogous portion of the proof in Theorem 1, as there is no need for
H\"older's inequality for integrals, and Lemma 5 can be used directly. The same bound is found to hold
for $\sum_{j\ge1}{\Cal I}_j$, with an additional log-factor.

The main contribution is a multiple of
$$
 \sum_{n\le N}d(n)n^{-3/4}\int_{2\pi\sqrt{nX}}^{2\pi\sqrt{2nX}}x_0^{1/2}\E^{iF(x_0)}|F''(x_0)|^{-1/2}\cdot
\frac{\z^4(\textstyle{\frac{1}{4}}+it)}{\frac14+it}\,\d t.\leqno(4.8)
$$
With $k=4, s = \frac14 +it$ Lemma 4 gives
$$
\eqalign{&
\z^4(s) = \chi^4(s)\z^4(1-s)\cr&
= \chi^4(s)\sum_{m=1}^\infty \rho(m/x)d_4(m)m^{s-1} + \sum_{m=1}^\infty \rho(m/y)d_4(m)m^{-s}\cr&
+ O(t^{1/3}) + O(t^{-2}y^{3/4}\log^3X).\cr}
$$
Here $xy = (t/(2\pi))^4$, and we choose $x = y = (t/(2\pi))^2$. The error terms are
$$
\ll t^{1/3} + t^{-1/2}\log^3X \ll t^{1/3},
$$
and the contribution to (4.8) is
$$
\ll X^{1/2}\sum_{n\le N}d(n)n^{-3/4}X^{3/4}n^{-1/4}(nX)^{1/6} \ll X^{17/12}N^{1/6}\log X.
$$
There remains the contribution of
$$
\sum_{n\le N}d(n)n^{-3/4}(I' + I''),
$$
where ($s = \frac14 + it$)
$$
\eqalign{
I' &:= \int_{2\pi\sqrt{nX}}^{2\pi\sqrt{2nX}}x_0^{1/2}\frac{\E^{iF(x_0)}}{\sqrt{|F''(x_0)|}}
\frac{\chi^4(s)}{s}\sum_{m\le 2x}\rho(m/x)d_4(m)m^{s-1}\,\d t\cr
I''&:= \int_{2\pi\sqrt{nX}}^{2\pi\sqrt{2nX}}x_0^{1/2}\frac{\E^{iF(x_0)}}{s\sqrt{|F''(x_0)|}}
\sum_{m\le 2x}\rho(m/x)d_4(m)m^{-s}\,\d t.\cr}
$$
Estimating, similarly as in the proof of Theorem 1, both $I'$ and $I''$ by the first derivative test
we obtain
$$
\eqalign{
I'&\ll X^{1/2}X^{3/4}n^{-1/4}(nX)^{1/4}\log^2X\ll X^{3/2}\log^2 X,\cr
I''&\ll X^{1/2}X^{3/4}n^{-1/4}(nX)^{-1/2}(nX)^{3/4}\log^2X\ll X^{3/2}\log^2 X.
\cr}
$$
This gives
$$
\sum_{n\le N}d(n)n^{-3/4}(I' + I'') \ll \sum_{n\le N}d(n)n^{-3/4}X^{3/2}\log^2X \ll X^{3/2}N^{1/4}\log^3X.
$$
Collecting all the estimates we infer that
$$
I(X) \ll_\e X^{13/8+\e}N^{-1/4} + X^{3/2}N^{1/4}\log^3X + X^{3/2+\e} +  X^{17/12}N^{1/6}\log X.
$$
We have $X^{3/2}N^{1/4} = X^{13/8}N^{-1/4}$ for $N = X^{1/4}$, which finally
gives
$$
I(X) \;\ll_\e\; X^{25/16+\e},
$$
as asserted in (1.7) of Theorem 2.

\head
5. Proof of Theorem 3
\endhead
Note that, by (1.1), one has
$$
\D(x) \;=\; \sum_{n\le x}d(n) - x(\log x + 2\gamma - 1),
\leqno(5.1)
$$
so that a comparison with (1.9) shows that the expression for $E(x)$ has a factor of $2\pi$ in the logarithm,
which (5.1) does not. It was observed first by M. Jutila \cite{Jut}, that the analogy between $\D(x)$ and
$E(x)$ becomes more exact if, instead of $\D(x)$, one introduces the function
$$
\D^*(x) := -\D(x)  + 2\D(2x) - \hf\D(4x)
= \hf\sum_{n\le4x}(-1)^nd(n) - x(\log x + 2\gamma - 1),
\leqno(5.2)
$$
Namely, then the function
$$
E^*(t):= E(t) - 2\pi \D^*\Bigl(\frac{t}{2\pi}\Bigr)\leqno(5.3)
$$
is smaller in mean square than either $\D^*(x)$ or $E(x)$. Namely Jutila \cite{Jut} proved that
$$
\int_T^{T+H}(E^*(t))^2\d t \,\ll_\e\, HT^{1/3}\log^3T + T^{1+\e}\qquad(2 \le H \le T), 
$$
while $E(x)$ satisfies an asymptotic formula analogous to (1.4). The first author \cite{Iv5} obtained the
asymptotic formula
$$
\int_0^T (E^*(t))^2\d t \;=\; T^{4/3}P_3(\log T) + O_\e(T^{7/6+\e}).\leqno(5.4)
$$
A comparison with (1.4) shows that in (5.4) the main term has the exponent smaller by 1/6.

Note that the analogue of Lemma 1 will hold with $\D^*(x)$ in place of $\D(x)$ (see e.g., (15.68)
of \cite{Iv1} for the analogue of (2.1)). Further, on using (5.4) we have
$$
\eqalign{&
\int_1^X E^*(x)\D_3(x)\d x  \ll \left\{\int_1^X (E^*(x))^2\d x\int_1^X \D_3^2(x)\d x\right\}^{1/2}\cr&
\ll X^{3/2}\log^{3/2}X.\cr}
$$
Here we have $13/9 < 3/2$, where 13/9 was the exponent in Theorem 1. Similarly it is found that
$$
\int_1^X E^*(x)\D_4(x)\d x \ll_\e X^{37/24+\e},
$$
but here $37/24 < 25/16$, where the latter is the exponent in our Theorem 2.

\medskip
Therefore, by (5.4), it remains to check that the integrals
$$
\int_X^{2X}\D^*\Bigl(\frac{x}{2\pi}\Bigr)\D_3(x)\d x,
\quad \int_X^{2X}\D^*\Bigl(\frac{x}{2\pi}\Bigr)\D_4(x)\d x
$$
satisfy the same bounds as the integrals in Theorem 1 and Theorem 2, respectively.
This is analogous as the  previous reasoning, with the only real difference that now the saddle point $x_0$
(we retain the same notation) will be different. Indeed, instead of  (3.6)
we shall have
$$
\eqalign{&
2\sum_{n\le N}(-1)^n d(n)n^{-3/4}\Re \Bigl\{\int_0^X\frac{\z^3(\frac14 +it)}{\frac14 +it}\times\cr&
\times \int_X^{2X}\,{\Bigl(\frac{x}{2\pi}\Bigr)}^{1/2+it}
\cos\Bigl(4\pi\sqrt{\frac{nx}{2\pi}} - \frac{\pi}{4}\Bigr)\d x\d t\Bigr\}.\cr}
$$
The factor $(-1)^n$ in the sum over $n$ is unimportant. Then we obtain, instead of $F(x)$,  the function
$$
\eqalign{
F^*(x) :&= \sqrt{8\pi nx} - t\log x,\cr
(F^*(x))'&= \sqrt{\frac{2\pi n}{x}} - \frac{t}{x},\cr
(F^*(x))'' &=\frac{t}{x^2} - \sqrt{\frac{\pi n}{2x^3}}.\cr}
$$
The new saddle point is
$$
x_0 = \frac{t^2}{2\pi n}.
$$
Then
$$
F^*(x_0) = 2t - 2t\log t + t\log(2\pi n).
$$
Instead of the functions $H_1, H_2$, we shall have now the functions $H_1^*, H_2^*$ satisfying
$$
\eqalign{
H_1^*(t) &= 3t\log\Bigl(\frac{2\pi\E}{t}\Bigr) + t\log m + F^*(x_0),\cr
(H_1^*(t))'&= -5\log t + \log(2\pi)^5mn,\cr
H_2^*(t)&= -t\log m + F^*(x_0),\cr
(H_2^*(t))'&= \log 2\pi\frac{n}{m} - 2\log t.\cr}
$$
In the relevant range for $m,n,t$ we then obviously have
$$
(H_1^*(t))' \,\gg\, \log t,\qquad (H_2^*(t))' \,\gg\, \log t,
$$
similarly as we had for $H'_1(t), H'_2(t)$.
Thus the proof goes through as before and finally one ends up with the bounds (1.10), as asserted.

\medskip
\vfill
\eject
\topglue1cm
\bigskip
\Refs

\smallskip

\item{[BW]}	J. Bourgain and N. Watt, Mean square of zeta function, circle problem
and divisor problem revisited, preprint available at {\tt arXiv:1709.04340}.

\smallskip
\item{[Iv1]} A. Ivi\'c, The Riemann zeta-function, John Wiley \&
Sons, New York 1985 (reissue,  Dover, Mineola, New York, 2003).

\smallskip
\item{[Iv2]} A. Ivi\'c, Large values of certain number-theoretic error terms, Acta Arith.
{\bf56}(1990), 135-159.

\smallskip
\item{[Iv3]} A. Ivi\'c,  Mean values of the Riemann zeta function, Tata Institute of Fund.
Research, LN's {\bf82}, (distributed by Springer Verlag, Berlin etc.), Bombay, 1991.
To be found online at
{\tt www.math.tifr.res.in/~publ/ln/tifr82.pdf}

\smallskip
\item{[Iv4]} A. Ivi\'c, Some problems on mean values of the Riemann zeta-function, Journal
    de Th\'eorie des Nombres Bordeaux {\bf8}(1996), 101-122.

\smallskip
\item{[Iv5]} A. Ivi\'c, On the mean square of the zeta-function and the divisor problem,
Annales Acad. Scien. Fennicae Math. {\bf32}(2007), 1-9.

\smallskip
\item{[Iv6]} A. Ivi\'c, On the divisor function and the Riemann zeta-function in short intervals,
The Ramanujan Journal  {\bf19}(2009), 207-224.

\smallskip
\item{[IvZh]} A. Ivi\'c and W. Zhai, On some mean value results for the zeta-function and
a divisor problem II, Indagationes Mathematicae, Volume {\bf26}, Issue 5,  2015,
pp. 842-866.

\smallskip
\item{[Jut]} M. Jutila, Riemann's zeta-function and the divisor problem,
Arkiv Mat. {\bf21}(1983), 75-96 and II, ibid. {\bf31}(1993), 61-70.

\smallskip
\item{[KV]} A.A. Karatsuba and S.M. Voronin, The Riemann zeta-function, Walter de Gruyter,
Berlin--New York, 1992.

\smallskip
\item{[Kol]} G. Kolesnik, On the estimation of multiple exponential sums. Recent progress in analytic number theory,
Vol. 1 (Durham, 1979), pp. 231--246, Academic Press, London-New York, 1981.

\smallskip
\item{[LT]}
Y.-K. Lau and K.-M. Tsang,  Mean square of the remainder term in the Dirichlet divisor problem,
 Journal de Th\'eorie des
Nombres de Bordeaux {\bf7}(1995), 75-92.

\smallskip
\item{[Meu]}
T. Meurman,
On the mean square of the Riemann zeta-function,
Quarterly J. Math., Oxford II. Ser. {\bf38}(1987), 337-343.

\smallskip
\item{[Tit]}  E.C. Titchmarsh, The theory of the Riemann
zeta-function (2nd edition),  Oxford University Press, Oxford, 1986.

\item{[Ton]} K.-C. Tong, On divisor problems III, Acta Math. Sinica {\bf6}(1956), 515-541.

\smallskip
\item{[Zh]}
W.-P. Zhang,
On the divisor problem, Kexue Tongbao (English Ed.) {\bf33}(1988), no. 17, 1484-1485.

\endRefs

\enddocument

\bye